\documentclass[12pt, reqno]{amsart}
\usepackage{amsmath, amsthm, amscd, amsfonts, amssymb, graphicx, color}
\usepackage[bookmarksnumbered, colorlinks, plainpages]{hyperref}
\usepackage{theoremref}
\hypersetup{colorlinks=true,linkcolor=red, anchorcolor=green, citecolor=cyan, urlcolor=red, filecolor=magenta, pdftoolbar=true}

\newtheorem{theorem}{Theorem}[section]
\newtheorem{lemma}[theorem]{Lemma}
\newtheorem{proposition}[theorem]{Proposition}
\newtheorem{corollary}[theorem]{Corollary}
\theoremstyle{definition}
\newtheorem{definition}[theorem]{Definition}

\newtheorem{remark}[theorem]{Remark}
\theoremstyle{remark}

\numberwithin{equation}{section}

\begin{document}
\setcounter{page}{1}

\title[]
{Characterization of positive definite, radial functions on free groups}

\author{Chian Yeong Chuah, Zhen-Chuan Liu, Tao Mei}

\address{Chian Yeong Chuah, Department of Mathematics, Ohio State University, 231 W. 18th Ave.
Columbus, OH 43210, USA.}
\email{\textcolor[rgb]{0.00,0.00,0.84}{chuah.21@osu.edu}}

\address{Zhen-Chuan Liu, Department of Mathematics, Baylor University, 1301 S University Parks Dr, Waco, TX 76798, USA.}
\email{\textcolor[rgb]{0.00,0.00,0.84}{zhen-chuan\_liu1@baylor.edu}}

\address{Tao Mei, Department of Mathematics, Baylor University, 1301 S University Parks Dr, Waco, TX 76798, USA.}
\email{\textcolor[rgb]{0.00,0.00,0.84}{Tao\_Mei@baylor.edu}}


\maketitle

\begin{abstract}
    This article studies the properties of positive definite, radial functions on free groups following the work of Haagerup and Knudby (\cite{HK2015}). We obtain characterizations of radial functions with respect to the $\ell^{2}$ length on the free groups with infinite generators and the characterization of the positive definite, radial functions with  respect to the $\ell^{p}$ length on the free real line with infinite generators for $0 < p \leq 2$. We obtain a L\'{e}vy-Khintchine formula for length-radial conditionally negative functions as well.
\end{abstract}

\section{Introduction}

Let $G$ be a group. A function $\varphi : G \to \mathbb{C}$ is called {\it positive definite} if the associated Toeplitz-type matrix $$[\varphi(x_i^{- 1}x_j)]_{1\leq i,j\leq n}$$ is positive definite for any $n\in{\Bbb N}$ and any  $(x_i)_{i=1}^n \in G$, i.e.
\begin{eqnarray*}
\sum_{i,j=1}^n\bar{c_i}c_j\varphi (x_i^{- 1}x_j)\geq0
\end{eqnarray*}  
for any complex numbers  $(c_i)_{i=1}^n$.
The classical  Bochner-Herglotz  theorem (\cite[5.5.2]{Bh07} says that a function on the integer group is positive definite if and only if it is the Fourier transform of a finite non-negative Borel measure on the torus.  

There is a similar concept of {\it positive definiteness} on semigroups. Let $G_+$ be a semigroup. A function $\dot\varphi : G_+ \to \mathbb{C}$ is called {\it positive definite} in the semigroup sense if the associated Hankel-type matrix $$[\dot{\varphi}(x_ix_j)]_{1\leq i,j\leq n}$$ is positive definite for any $n\in{\Bbb N}$ and any $n$ elements $x_i \in G_+$. This is equivalent to saying that $\sum_{i,j=1}^n\bar{c_i}c_j\dot{\varphi} (x_ix_j)\geq0$  for any complex numbers $c_i, 1\leq i\leq n$.  The Hamburger theorem (\cite[Chapter 1,Theorem 7.1]{Pe03} says that a bounded function $\dot\varphi$ is positive definite on the semigroup ${\Bbb Z}_+={\Bbb N}\cup\{0\}$ if and only if  $\dot\varphi$ is the moment of a nonnegative Borel measure $\mu$ on $[-1,1]$, i.e. there exists $\mu$ such that
\begin{eqnarray}\label{moment}
\dot\varphi(k)=\int_{-1}^1 t^k d\mu(t) .
\end{eqnarray}
for $k\in \mathbb{Z}_{+}$. Note that the support of $\mu$ is a subset of $[-1,1]$ as $\dot\varphi$ is bounded. Given any   bounded positive definite function $\dot\varphi$ on ${\Bbb Z}_+$, the formula $\varphi(k)=\dot\varphi(|k|)$ defines a symmetric positive definite function on ${\Bbb Z}$.  This can be seen by (\ref{moment}) and the well-known fact that $k\mapsto t^{|k|}$ is positive definite on ${\Bbb Z}$ for any $-1\leq t\leq 1$. However, not every  symmetric positive definite function $\varphi$ on ${\Bbb Z}$ is of the form of $\varphi(k)=\dot\varphi(|k|)$.  In fact, Haagerup and Knudby proved   in \cite[Theorem 3.3] {HK2015} that  there is a one to one correspondence between the class of bounded positive definite functions on ${\Bbb Z}_+$ and the class of radial positive definite functions on the infinite free product of ${\Bbb Z}$. The notion of radial (or spherical) functions is first introduced and studied in \cite{F-TP1982}. We restate Haagerup's result as follows.
\bigskip




\begin{theorem}[Haagerup-Knudby, \cite{HK2015}] Let $\mathbb{F}_{\infty}$ be the free group of countable many infinite generators.  Let $\| g \|_{1}$ be the reduced word length of an element $g\in{\Bbb F}_\infty$.  Given a bounded function $\dot\varphi$ on ${\Bbb Z}_+$, the following are equivalent.

\begin{enumerate}
\item The function $\varphi(g)=\dot{\varphi}(\| g \|_{1})$ is positive definite on ${\Bbb F}_\infty$.

\item There is a finite positive Borel measure $\mu$ on $[- 1, 1]$ such that 

\begin{displaymath}
\dot{\varphi}(n) = \int_{- 1}^{1} s^{n} \ d \mu (s), \ n \in \mathbb{N}
\end{displaymath}

\end{enumerate}
\end{theorem}
Together with the Hamberger theorem (\ref{moment}), Haagerup-Knudby's theorem gives a one to one correspondence between the class of bounded positive definite functions on ${\Bbb Z}_+$ and the class of radial positive definite functions on the infinite free product of ${\Bbb Z}$, though the article does not provide a direct proof of this correspondence.

In this article, we give a direct argument for this correspondence. Our argument  works  for  more general settings (see \thref{exact} ) including the $\ell_{p}$-length radial functions on the infinitely generated free real line for all $0<p\leq2$. The $p=2$ case can be viewed as a free analogue of the classical Schoenberg-Bochner theorem (see e.g. \cite[Theorem 13.14]{SSV2012})  which says that a function $f$ is conditionally negative definite on $[0, \infty)$
 in the semigroup sense  if and only if  the function $\xi \mapsto f(\| \xi \|^{2})$, $\xi \in \mathbb{R}^{d}$, is conditionally negative definite for all $d \in \mathbb{N}$.







\begin{definition} Fix a set of generators $\{g_i, i\in{\Bbb N}\}$ of the free group $ \mathbb{F}_{\infty}$. Let  $0 < p \leq 2$. For a reduced word $g= g_{i_1}^{k_{1}} g_{i_2}^{k_{2}} ... g_{i_n}^{k_{n}} \in \mathbb{F}_{\infty}$, define the $\ell^{p}$ length of $g$, denoted by $\left\| g \right \|_{p}$ as:

\begin{displaymath}
\left\| g \right \|_{p} = \left( \sum_{j = 1}^{n} | k_{j} |^{p} \right)^{\frac{1}{p}}
\end{displaymath}

\end{definition}

The maps $g\mapsto \left\| g \right \|_{p}^{p} $ are still conditionally negative definite (see \thref{cn}) on the free group  (see Proposition \ref{generallength}). We say a  function $\varphi$  on $\mathbb{F}_\infty$ is $\|\cdot\|_p$-radial if  the value of $\varphi(g)$ only depends on $\|g\|_p$. Our first main result is stated as follows.

\begin{theorem}
 Suppose $\varphi$ is a $\|\cdot\|_2$-radial function on $\mathbb{F}_\infty$  with $\varphi(e) = 1$, where $e$ is the identity of $\mathbb{F}_{\infty}$. The following are equivalent. 

\begin{enumerate}
\item $\varphi(g) $ defines a positive definite function on $\mathbb{F}_{\infty}$.

\item There is a probability measure $\mu$ on $[- 1, 1]$ such that 

\begin{displaymath}
\varphi(g) = \int_{-1}^{1} s^{\|g\|_{2}^{2}} \ d \mu(s)
\end{displaymath}

\end{enumerate}

Moreover, if (2) holds, then $\mu$ is uniquely determined by $\varphi$.
\end{theorem}

\begin{theorem} Suppose $\psi : \mathbb{F}_{\infty} \to \mathbb{C}$ is an $\| \cdot \|_{2}$-radial function with $\psi (e) = 0$, where $e$ is the identity of $\mathbb{F}_{\infty}$. Then, the following are equivalent: 

\begin{enumerate}
    \item $\psi$ is conditionally negative definite on $\mathbb{F}_{\infty}$
    
    \item There is a probability measure $\nu$ on $[- 1, 1]$ such that 
    
    \begin{displaymath}
    \psi (g) = \int_{- 1}^{1} \frac{1 - s^{\left\| g \right\|_{2}^{2}}}{1 - s} \ d \nu(s)
    \end{displaymath}
\end{enumerate}
 
Moreover, if (2) holds, then $\nu$ is uniquely determined by $\psi$. 

\end{theorem}









We obtain characterizations for $\|\cdot\|_p$-radial, positive definite and $\|\cdot\|_p$-radial conditionally negative definite functions on the infinite free product of the group of real numbers as well. See \thref{positive equivalence real} and \thref{negative equivalence real}. Similar results for the commutative case is also generalized in \thref{commutative result}.

Positive definite functions are closely connected to completely positive   maps of the Fourier multiplier type. Let $G$ be a group and $\varphi : G \to \mathbb{C}$ a bounded function. Let $\lambda_{s}$ be the left regular representation of $s \in G$. Consider the associated multiplier $M_{\varphi}$   on $Span[ \lambda(G) ]$  defined as:

\begin{equation}\label{multiplier}
    M_{\varphi} \left( \sum  c_{s} \lambda_{s} \right) =  \sum  \varphi(s) c_{s} \lambda_{s}
\end{equation}
Then $M_\varphi$ extends to a  completely positive map on $C_{r}^{*} (G)$ if and only if $\varphi$ is positive definite. In this case, $M_\varphi$ is also completely bounded on $C_{r}^{*} (G)$ with norm $\varphi(e)$. For more information on the complete positivity of the multipliers, the readers can refer to \cite{MS2017}.

Following Haagerup's pioneer work (\cite{Ha79}), the completely positivity and the completely boundedness of the map $M_\varphi$, with $\varphi$ being a radial function with respect to the $\ell_1$-length, are fully characterized (\cite{HSS2010, HK2015}). For more information about the study of positive definite functions on free groups, the reader can refer to (\cite{B1989} and \cite{B1987}). These works significantly improve  the understanding of the approximation properties of the free groups and the associated noncommutative $L^p$-spaces. Nevertheless, our understanding is still incomplete. For instance, the existence of a Schauder basis for the reduced free group  $C^*$ algebra is still a mystery. Bo\.zejko and Fendler's work in (\cite{BF2006}) implies that, the sequence of $\{ l_{g} : g \in \mathbb{F}_{n} \}$, enumerated in an order compatible with the word length, is not a Schauder basis of the non-commutative $L^{p}$ spaces associated with the von Neumann algebra of the rank $n$ free group $\mathbb{F}_{n}$ if $p > 3$ or $1\leq p<\frac{3}{2}$ and $n \geq 2$. A better understanding of positive definite functions beyond $\| \cdot \|_{1}$-radial type would help.  The main results of this article (Theorem 1.2 and 1.3)  complement Haagerup and Knudby's work (\cite{HK2015}) and provide characterizations of the complete positivity of the corresponding multipliers $M_\varphi$ defined as in (\ref{multiplier}) with $\varphi$ being $\| \cdot \|_{2}$-radial. The classical $\| \cdot \|_{2}$-radial Fourier multipliers are those associated with the Laplacian operators. We hope the results obtained in this article will shed light on determining an appropriate  Laplace type operator on the free group $C^*$-algebras.

\section{Preliminaries}

First, we recall the general definition of a positive definite function and a conditionally negative definite function.

\begin{definition} Let $G$ be a group. A function $\varphi: G \to \mathbb{C}$ is Hermitian if $\varphi(g^{- 1}) = \overline{\varphi(g)}$ for all $g \in G$.  

\end{definition}


\begin{definition} Let $G$ be a group. A function $\varphi : G \to \mathbb{C}$ is positive definite if for each $n \in \mathbb{N}$, $\{ x_{1}, ..., x_{n} \} \subseteq G$ and $\{ c_{1}, ..., c_{n} \} \subseteq \mathbb{C}$, 

\begin{displaymath}
\sum_{j, k = 1}^{n} c_{j} \ \overline{c_{k}} \ \varphi(x_{j}^{- 1} x_{k}) \geq 0.
\end{displaymath}

\end{definition}

\begin{definition} \thlabel{cn} Let $G$ be a group. A function $\psi : G \to \mathbb{C}$ is conditionally negative definite if 

\begin{enumerate}
    \item $\psi$ is Hermitian.
    
    \item For each $n \in \mathbb{N}$, $\{ x_{1}, ..., x_{n} \} \subseteq G$ and $\{ c_{1}, ..., c_{n} \} \subseteq \mathbb{C}$ satisfying 

\begin{displaymath}
\sum_{j = 1}^{n} c_{j} = 0, \text{ we have that } \sum_{j, k = 1}^{n} c_{j} \ \overline{c_{k}} \ \varphi(x_{j}^{- 1} x_{k}) \leq 0.
\end{displaymath}
\end{enumerate} 

\end{definition}

\begin{lemma}[Bochner's Theorem,see \cite{WR1990},page 19]
 Let $G$ be a locally compact, abelian group and let $\phi : G \to \mathbb{C}$ be a continuous function. Then, $\phi$ is positive definite on $G$ if and only if there exists a non-negative, finite Radon measure $\mu$ on the dual group $\Gamma$ of $G$ such that

\begin{equation*}
    \phi(x) = \int_{\Gamma} \gamma(x) \ d \mu(\gamma) \ \ \ \ \ \  (x \in G),
\end{equation*}
where $\gamma$ is the character function on $G$.
\end{lemma}
Now, we recall Schoenberg's theorem which characterizes the conditionally negative definite functions on $G$, see \cite{BCR1984}.
\begin{lemma}[Schoenberg]
 Let $G$ be a group. Let $\psi : G \to \mathbb{C}$ be a Hermitian function. Then, the following are equivalent.

\begin{enumerate}
\item $\psi$ is conditionally negative definite on $G$.

\item For each $t > 0$, the function $\varphi_{t} : G \to \mathbb{C}$ defined by $\varphi_{t} (g) := e^{- t \psi(g)}$ is positive definite. 
\end{enumerate}
\end{lemma}
Besides Schoenberg's theorem, there is also another classical result which relates a conditionally negative definite kernel with a function on a Hilbert space. However, this requires some additional assumption that the kernel is real valued and is zero on the diagonal.  

\begin{lemma} Let $G$ be a group. Let $\psi : G \to \mathbb{R}$ be a real-valued function where $\psi (e) = 0$. Then, the following are equivalent.

\begin{enumerate}
\item $\psi$ is conditionally negative definite on $G$.

\item There exist a Hilbert space $H$ and a function $f : G \to H$ such that $\psi (x^{- 1} y) = \left\| f(y) - f(x) \right\|^{2}$ for all $x$, $y \in G$. 
\end{enumerate}

\end{lemma}



\begin{definition} Let $G$ be a group and $\theta : G \to \mathbb{R}_{+}$ be a function. A function $\varphi : G \to \mathbb{C}$ is said to be radial with respect to $\theta$ if there exists $\dot{\varphi} : \text{ran}(\theta) \to \mathbb{R}_{+}$ such that for all $g \in G$, $\varphi (g) = \dot{\varphi}[\theta (g)]$.
\end{definition}



Apart from the case for groups, there is also an analogous definition of positive definite functions in the setting of an abelian semi-group. However, care must be taken since in general, a semi-group does not have an inverse. As a remark, one can prove the theory for a general involution semigroup. However, for our purpose, we always assume the involution operator to be the identity operator.

\begin{definition} Let $S$ be an abelian semigroup. A function $\varphi : S \to \mathbb{C}$ is positive definite in the semi-group sense if for each $n \in \mathbb{N}$, $\{ s_{1}, ..., s_{n} \} \subseteq \mathbb{N}$ and $\{ c_{1}, ..., c_{n} \} \subseteq \mathbb{C}$, 

\begin{equation*}
\sum_{j, k = 1}^{n} c_{j} \ \overline{c_{k}} \ \varphi( s_{j} + s_{k}) \geq 0.
\end{equation*}

\end{definition}

\begin{definition} Let $S$ be an abelian semigroup. A function $\psi : S \to \mathbb{C}$ is conditionally negative definite in the semi-group sense if 

\begin{enumerate}
    \item $\psi$ is real-valued.
    
    \item For each $n \in \mathbb{N}$, $\{ s_{1}, ..., s_{n} \} \subseteq S$ and $\{ c_{1}, ..., c_{n} \} \subseteq \mathbb{C}$ satisfying 

\begin{displaymath}
\sum_{j = 1}^{n} c_{j} = 0, \text{ we have that } \sum_{j, k = 1}^{n} c_{j} \ \overline{c_{k}} \ \varphi(s_{j} + s_{k}) \leq 0.
\end{displaymath}
\end{enumerate} 

\end{definition}






\begin{definition} A function $\rho : S \to \mathbb{R}$ is called a semicharacter on the  abelian semigroup $S$ if

\begin{enumerate}
    \item $\rho (0) = 1$.
    \item $\rho (s + t) = \rho(s) \rho (t)$ for $s$, $t \in S$.
\end{enumerate}
\end{definition}

\begin{definition} Let $S$ be an abelian semigroup. The set $S^{*} := \{ \rho : \rho$ is a semicharacter$\}$ equipped with the topology of pointwise convergence is called the dual semigroup of $S$.
\end{definition}
It turns out that being equipped with the topology of pointwise convergence, $S^{*}$ becomes a completely regular space, in particular it is a Hausdorff space. Moreover, it forms a topological semigroup, with the multiplication defined via pointwise multiplication and the constant function $1$ being the identity.  

\begin{definition}: Let $S$ be an abelian semigroup. The set $\widehat{S} := \{ \rho \in S^{*} : | \rho (s) | \leq 1 $ for $s \in S \}$ is called the restricted dual semigroup.
\end{definition}
By inheriting the subspace topology of $S^{*}$, $\widehat{S}$ becomes a compact subsemigroup of $S^{*}$, we need the following result.
\begin{lemma}[\cite{BCR1984},Theorem 4.2.8, page 96]\thlabel{semi} Let $S$ be an abelian semigroup. A function  $\varphi : S \to \mathbb{C}$ is positive definite and bounded on $S$ if and only if
\begin{equation*}
    \varphi (s) = \int_{\widehat{S}} \rho (s) \ d \mu (\rho) \ \ \ \ \ \ \ (s \in S),
\end{equation*}
where $\mu$ is a Radon measure on $\widehat{S}$. Moreover, if we assume that $\varphi (e) = 1$, then $\mu$ is a probability measure. 
\end{lemma}

\begin{remark}
Note that there is a slight difference between the formulation of the results for case of groups and semigroups.
\end{remark}

\section{Key Lemma}



\begin{definition} Let $G$ be a group and $\theta : G \to \mathbb{R}_{+}$ be a function on $G$ such that the range of $\theta$, $ran(\theta)$, is a semi-group. Then, $\theta$ is said to be a partial morphism on $G$ if it satisfies the following property:

Given $N \in \mathbb{N}$, $M \in \mathbb{N}$ and $s_{1} , \cdots , s_{M} \in ran(\theta)$, there exist elements $\{ g_{n, s_{k}} \in G : 1 \leq n \leq N \text{ and } 1 \leq k \leq M \}$ such that $\theta(g_{n, s_{j}}^{- 1} g_{m, s_{k}}) = s_{k} + s_{j}$ for $n \not = m$.

\end{definition}

\begin{lemma}\thlabel{exact} Let $\theta : G \to \mathbb{R}_{+}$ be a partial morphism. Let $\varphi : G \to \mathbb{R}$ be a positive definite, radial function with respect to $\theta$.  Then, the corresponding function $\dot{\varphi} : ran(\theta) \to \mathbb{R}$ satisfying $\dot{\varphi} [\theta(g)] = \varphi (g)$ is positive definite in the semi-group sense and is bounded.
\end{lemma}
\begin{proof}
Fix $N \in \mathbb{N}$, $K \in \mathbb{N}$. Consider $\{ s_{1}, \cdots, s_{K}\} \subseteq ran(\theta)$, $\{ c_{1}, \cdots, c_{K} \} \subseteq \mathbb{C}$. 

Next, we set $d_{1, 1} = d_{2, 1} = \cdots = d_{N, 1} :\equiv c_{1}$, $d_{1, 2} = d_{2, 2} = \cdots = d_{N, 2} :\equiv c_{2}$, $\cdots, d_{1, K} = d_{2, K} = \cdots = d_{N, K} :\equiv c_{K}$.

Since $\theta$ is a  partial morphism on $G$, there exist  elements $\{ g_{n, s_{k}} \in G : 1 \leq n \leq N \text{ and } 1 \leq k \leq K \}$ such that $\theta(g_{n, s_{j}}^{- 1} g_{m, s_{k}}) = s_{k} + s_{j}$ for $n \not = m$.

Since $\varphi$ is positive definite on $G$, $\sum_{n, m = 1}^{N} \sum_{k, j= 1}^{K} d_{m, k} \overline{d_{n, j}} \varphi (g_{n, s_{j}}^{- 1} g_{m, s_{k} } ) \geq 0$. Next, we perform some calculation. 

For each $1 \leq n, m \leq N$, we have that
\[
    \sum_{j, k = 1}^{K} d_{m, k} \overline{d_{n, j}} \varphi (g_{n, s_{j}}^{- 1} g_{m, s_{k}}) = \sum_{j, k = 1}^{K} d_{m, k} \overline{d_{n, j}} \dot{\varphi} [\theta (g_{n, s_{j}}^{- 1} g_{m, s_{k}} )].
\]
If $m \not = n$, then
\[
    \sum_{j, k = 1}^{K} d_{m, k} \overline{d_{n, j}} \varphi (g_{n, s_{j}}^{- 1} g_{m, s_{k}}) = \sum_{j, k = 1}^{K} d_{m, k} \overline{d_{n, j}} \dot{\varphi} [\theta (g_{n, s_{j}}^{- 1} g_{m, s_{k}} )] = \sum_{j, k = 1}^{K} d_{m, k} \overline{d_{n, j}} \dot{\varphi} (s_{j} + s_{k}).
\]
Taking the sum of $m$, $n$ from $1$ to $N$, we obtain
\begin{align*}
&\sum_{m, n = 1}^{N} \sum_{j, k = 1}^{K} d_{m, k} \overline{d_{n, j}} \varphi (g_{n, s_{j}}^{- 1} g_{m, s_{k}}) \\
&=  \sum_{n = 1}^{N} \sum_{j, k = 1}^{K} d_{n, k} \overline{d_{n, j}} \varphi (g_{n, s_{j}}^{- 1} g_{m, s_{k}}) +  \sum_{1 \leq m \not = n \leq N } \sum_{k, j = 1}^{K} d_{m, k} \overline{d_{n, j}} \dot{\varphi} (s_{j} + s_{k}).
\end{align*}

\noindent Using the relationship between $d_{n, k}$ and $c_{k}$, we have
\begin{align*}
0 \leq & \sum_{m, n = 1}^{N} \sum_{j, k = 1}^{K} d_{m, k} \overline{d_{n, j}} \varphi (g_{n, s_{j}}^{- 1} g_{m, s_{k}}) = \sum_{m, n = 1}^{N} \sum_{j, k = 1}^{K} c_{k} \overline{c_{j}} \varphi (g_{n, s_{j}}^{- 1} g_{m, s_{k}})\\
= &\sum_{n = 1}^{N} \sum_{j, k = 1}^{K} c_{k} \overline{c_{j}} \varphi (g_{n, s_{j}}^{- 1} g_{n, s_{k}}) + (N^{2} - N) \sum_{k,j=1}^{M} c_{k} \overline{c_{j}} \dot{\varphi} (s_{j} + s_{k}) \\
\leq & \sum_{n = 1}^{N} \sum_{j, k = 1}^{K} |c_{k}| |c_{j}| |\varphi(g_{n, s_{j}}^{- 1} g_{n, s_{k}})| + (N^{2} - N) \sum_{k,j=1}^{M} c_{k} \overline{c_{j}} \dot{\varphi} (s_{j} + s_{k}) \\  
\leq & \sum_{n = 1}^{N} \sum_{j, k = 1}^{K} |c_{k}| |c_{j}| \varphi(0) + (N^{2} - N) \sum_{k,j=1}^{M} c_{k} \overline{c_{j}} \dot{\varphi} (s_{j} + s_{k}) \\
=& N \sum_{j, k = 1}^{K} |c_{k}| |c_{j}| \varphi(0) + N (N - 1) \sum_{k,j=1}^{M} c_{k} \overline{c_{j}} \dot{\varphi} (s_{j} + s_{k}) .
\end{align*}
\text{Thus, }
\begin{align*}
 - \frac{1}{N - 1} \left[ \sum_{j, k = 1}^{K} |c_{k}| |c_{j}| \varphi(0) \right] \leq \sum_{k, j = 1}^{M} c_{k} \overline{c_{j}} \dot{\varphi} (s_{j} + s_{k}).
\end{align*}

Since the above inequality holds true for all $N \in \mathbb{N}$, we have $\sum_{k,j=1}^{M} c_{k} \overline{c}_{j} \dot{\varphi}(s_{j} + s_{k} ) \geq 0$. Since the above inequality holds true for any $\{ s_{1}, \cdots, s_{K}\} \subseteq G$ and any $\{ c_{1}, \cdots, c_{K} \} \subseteq \mathbb{C}$,  $\dot{\varphi}$ is positive definite in the semi-group sense on $S$.
Finally, $\dot{\varphi}$ is bounded because $\varphi$ is bounded. 

\end{proof}

\begin{remark}
It turns out that \thref{exact} still holds true under a slightly weaker assumption. We will state the theorem and leave the proof to the readers. 
\end{remark}

\begin{lemma}
Given a group $G$, suppose the function $\theta : G \to \mathbb{R}_{+}$ satisfies that for any subset 
 $ \ U = \{ s_{1}, \cdots, s_{K} \}  $ of ${ran}(\theta)$,  there exists  a sequence $$V_{N} = \{ g_{n, s_{k}} \in G : 1 \leq n \leq N, 1 \leq k \leq K \}, N\in{\Bbb N}$$ such that 

\begin{eqnarray*}
\lim_{N \to \infty} \frac{ \#  \{ (m, n ) \in \mathbb{N}^{2} : 1 \leq m, n \leq N, \theta (g_{m, j}^{- 1} g_{n, k}) = s_{j} + s_{k},    \forall  1 \leq j, k \leq K \}}{N^{2}}  = 1.
\end{eqnarray*}

Then, for each positive definite function $\varphi : G \to \mathbb{R}$ which is radial with respect to $\theta$, the corresponding function $\dot{\varphi} : \text{ran} (\theta) \to \mathbb{R}$ satisfying $\dot{\varphi} [\theta(g)] = \varphi(g)$ is positive definite in the semi-group sense. 
\end{lemma}

\begin{corollary} \thlabel{exact negative definite} Let $\theta : G \to \mathbb{C}$ be a function. Let $\psi : G \to \mathbb{R}$ be a conditionally negative definite, radial function with respect to $\theta$. Suppose that $\theta$ is a partial morphism on $G$. Then, the corresponding function $\dot{\psi} : {\rm ran}( \theta) \to \mathbb{R}$ satisfying $\dot{\psi} [\theta(g)] = \psi (g)$ is conditionally negative definite in the semi-group sense and bounded below.
\end{corollary}
\begin{proof} Since $\psi$ is a conditionally negative definite function on the group $G$, by Schoenberg's theorem, for all $t > 0$, the function $\varphi_{t} : G \to \mathbb{C}$ defined by $\varphi_{t} (g) := e^{- t \psi (g)}$ is positive definite on $G$. By \thref{exact}, the function $\dot{\varphi_{t}} : ran(\theta) \to \mathbb{R}$ defined by $\dot{\varphi_{t}} [\theta(g)] := e^{- t \dot{\psi}[\theta(g)]}$ is positive definite on $ran(\theta)$. Again, by Schoenberg's theorem, the function $\dot{\psi}$ is conditionally negative definite in the semi-group sense on ${\rm ran}({\theta})$.
Finally, $\dot{\psi}$ is bounded below because $\psi$ is bounded below.
\end{proof}

\begin{theorem}\thlabel{main} Let $G$ be a group and $\theta : G \to \mathbb{R}_{+}$ be a partial morphism on $G$. Let $\varphi : G \to \mathbb{R}$ be a radial function with respect to $\theta$ where $\varphi(e) = 1$, If $\varphi$ is positive definite on $G$, then there exists a unique probability measure $\mu$ on $\widehat{ran(\theta)}$ such that 

\begin{equation*}
\varphi(g)=\dot{\varphi} [ \theta(g) ]=\int_{\widehat{S} } \rho [\theta(g)] \ d\mu(\rho).
\end{equation*}

\end{theorem}
\begin{proof} Since $\varphi$ is positive definite on $G$, $\dot{\varphi}$ is positive definite in the semi-group sense on $ran(\theta)$ by \thref{exact}. Then, by \thref{semi}, there exists a unique probability measure $\mu$ on $\widehat{S}$ such that 

\begin{equation*}
\varphi(g)=\dot{\varphi} [ \theta(g) ] =\int_{\widehat{S} } \rho [\theta(g)] \ d\mu(\rho).
\end{equation*} 
\end{proof}

\section{Proof of Main Theorems}

\subsection{Case of $\ell^{2}$ length for $\mathbb{F}_{\infty}$} \ 

We consider the case where the group $G = \mathbb{F}_{\infty}$, the free group with infinite generators, and the function $\theta = \left\| \cdot \right\|_{2}^{2}$ is the $\ell^{2}$ length of an element in $\mathbb{F}_{\infty}$.

In this case, $ran(\theta) = \mathbb{N}$ and hence, $\widehat{ran(\theta)} \cong [ - 1, 1 ]$. More precisely, given $\rho \in \widehat{ \mathbb{N}}$, there exists a unique $x \in [ - 1, 1 ]$ such that $\rho (n) = x^{n}$. First, we show that the $\ell^{2}$ length function is indeed a partial morphism on $\mathbb{F}_{\infty}$. 

\begin{proposition} The function $\| \cdot \|_{2}^{2} : \mathbb{F}_{\infty} \to \mathbb{N}$ is a partial morphism on $\mathbb{F}_{\infty}$.
\end{proposition}
\begin{proof}
 We label the infinite generators of $\mathbb{F}_{\infty}$ as the entries of the following infinite matrix 
 \[
\begin{pmatrix}
 g_{1,1} & \cdots & g_{1,n} & \cdots \\
 \vdots  & \ddots &  \vdots & \cdots \\
 g_{n,1} &  \cdots  & g_{n,n} & \cdots \\
 \vdots  &        & \vdots  & \ddots
\end{pmatrix}.
\]
For each $n \in \mathbb{N}$, define $q_n: \mathbb{N}^+ \rightarrow \mathbb{F}_\infty$ by
$$
q_n(j) := g_{n, j}g_{n, j-1} \cdots g_{n, 1}.
$$
Then, for each $m$, $n$, $j$, $k \in \mathbb{N}$, there is
\begin{align*}
\|[q_n(j)]^{-1}q_m(k)\|_{2}^{2}
=
\begin{cases}
0,
& \text{if } ~m=n~and~k=j,\\
j+k, & \text{otherwise}. 
\end{cases} 
\end{align*}
Hence, the proposition is proved.
\end{proof}

\begin{proposition}\thlabel{length}
Let $\mathbb{F}_r$ be a free group with generators $g_1, g_2, \cdots, g_{r}$, where $r \in \mathbb{N} \cup \{ \infty \}$.
\begin{enumerate}
\item  Let $s\in[0, 1]$. Then, the function $\psi : \mathbb{F}_{r} \to \mathbb{R}$ defined by $\psi(g) := s^{\|g\|_{2}^{2}}$ is positive definite on $\mathbb{F}_r$, i.e. the function $\varphi : \mathbb{F}_{r} \to \mathbb{R}$ defined by $\varphi (g) := \|g\|_2^2$ is conditionally negative definite on $\mathbb{F}_r$.

\item The function $\varphi : \mathbb{F}_{r} \to \mathbb{R}$ defined by $\varphi(g) := (-1)^{\|g\|_{2}^{2}}$ is positive definite on $\mathbb{F}_r$.
\end{enumerate}
Moreover, $g\rightarrow s^{\|g\|_{2}^{2}}, s\in [-1, 1]$, is a positive definite function on $\mathbb{F}_r$.

\end{proposition}
\begin{proof}
 (1) follows from Proposition \ref{generallength}. (2) can be directly verified by definition.
Since a product of two positive definite functions is also positive definite, we conclude that for $- 1 \leq s \leq 1$, the function $g \rightarrow s^{\|g\|_{2}^{2}}$ is positive definite on $\mathbb{F}_r$.
\end{proof}









\begin{corollary} \thlabel{main positive definite}Given $\varphi : \mathbb{F}_{\infty} \to \mathbb{R}$ a radial function where $\varphi(e) = 1$, the following are equivalent. 

\begin{enumerate}
\item $\varphi$ is positive definite on $\mathbb{F}_{\infty}$.

\item There is a probability measure $\mu$ on $[- 1, 1]$ such that 

\begin{displaymath}
\varphi(g)=\dot{\varphi}(\|g\|_{2}^{2})= \int_{-1}^{1} s^{\|g\|_{2}^{2}} \ d \mu(s)
\end{displaymath}

\end{enumerate}

Moreover, if (2) holds, then $\mu$ is uniquely determined by $\varphi$.

\end{corollary}

\begin{proof} (1) $\Longrightarrow (2)$ follows from \thref{main}. To prove (2) $\Longrightarrow$ (1), let $\mu$ be a probability measure on $[- 1, 1]$. For each $s \in [- 1, 1]$, the function $\psi_{s} (g) := s^{\| g \|_{2}^{2}}$ is positive definite on $\mathbb{F}_{\infty}$ by \thref{length}. Taking finite sums and limits, we deduce that the function $\varphi : \mathbb{F}_{\infty} \to \mathbb{R}$ defined by $\varphi(g) := \int_{- 1}^{1} s^{\left\| g \right\|_{2}^{2}} \ d \mu(s)$ is positive definite. \end{proof}

\begin{theorem} Let $\psi : \mathbb{F}_{\infty} \to \mathbb{C}$ be an $\ell^{2}$-radial function where $\psi (e) = 0$. Then, the following are equivalent: 

\begin{enumerate}
    \item $\psi$ is conditionally negative definite on $\mathbb{F}_{\infty}$.
    
    \item There is a probability measure $\nu$ on $[- 1, 1]$ such that 
    
    \begin{displaymath}
    \psi (g) = \int_{- 1}^{1} \frac{1 - s^{\left\| g \right\|_{2}^{2}}}{1 - s} \ d \nu(s).
    \end{displaymath}
\end{enumerate}
 
Moreover, if (2) holds, then $\nu$ is uniquely determined by $\psi$. 

\end{theorem}
\begin{proof} $((1)\Longrightarrow (2))$
By Schoenberg's theorem, since $\psi$ is conditionally negative definite and $\psi(e) = 0$, the function $\varphi_{t} : \mathbb{F}_{\infty} \to \mathbb{C}$ defined by $\varphi_{t} (g) := e^{- t \psi (g)}$ is positive definite for each $t > 0$. Also, $\varphi_{t} (g) = 1$. By \thref{main positive definite}, there exists a unique probability measure $\mu_{t}$ on $[- 1, 1]$ such that

\begin{displaymath}
e^{- t \psi(g)} = \varphi_{t} (g) = \int_{- 1}^{1} s^{ \left\| g \right\|_{2}^{2}} \ d \mu_{t} (s).
\end{displaymath}

Now, let $t > 0$ and define a new measure on the Borel $\sigma$-algebra of $[- 1, 1]$, $\nu_{t}$ by $\nu_{t} (E) := \int_{- 1}^{1} \chi_{E}(s)  \frac{1 - s}{t} \ d \mu_{t}(s)$. Note that

\begin{align*}
&\frac{1 - e^{- t \psi(g)}}{t} = \int_{[- 1, 1)} \frac{1 - s^{ \left\| g \right\|_{2}^{2}}}{t} \ d \mu_{t}(s) + \int_{\{ 1 \}} \frac{1 - s^{ \left\| g \right\|_{2}^{2}}}{t} \ d \mu_{t}(s) \\
&= \int_{[- 1, 1)} \frac{1 - s^{ \left\| g \right\|_{2}^{2}}}{t} \ d \mu_{t}(s)  \quad ( \text{ since } \frac{1 - s^{ \left\| g \right\|_{2}^{2}}}{t} = 0 \text{ at } x = 1  )\\
&= \int_{[- 1, 1)} \frac{1 - s^{\left\| g \right\|_{2}^{2}}}{1 - s} \frac{1 - s}{t} \ d \mu_{t}(s) = \int_{[- 1, 1)}  \frac{1 - s^{\left\| g \right\|_{2}^{2}}}{1 - s} \ d \nu_{t}(s) \\
&= \int_{[- 1, 1)}  \frac{1 - s^{\left\| g \right\|_{2}^{2}}}{1 - s} \ d \nu_{t}(s) + \int_{\{ 1 \}}  \frac{1 - s^{\left\| g \right\|_{2}^{2}}}{1 - s} \ d \nu_{t}(s) \left( \text{ since } \nu_{t}(\{ 1 \}) = 0 \right) \\
&= \int_{[- 1, 1]}  \frac{1 - s^{\left\| g \right\|_{2}^{2}}}{1 - s} \ d \nu_{t}(s).
\end{align*}

Applying the identity for $h \in \mathbb{F}_{\infty}$ where $\left\| h \right\|_{2}^{2} = 1$, we obtain:

\begin{displaymath}
\nu_{t} ([- 1, 1]) = \frac{1 - e^{- t \psi(h)}}{t}
\end{displaymath}

Taking the supremum over all $t > 0$,
\begin{equation*}
    \sup_{t > 0} \left\| \nu_{t} \right\|_{\text{var}} =  \sup_{t > 0} \frac{1 - e^{- t \psi(h)}}{t} = \psi(h).
\end{equation*}

So, the set $\{ \nu_{t} : t > 0 \}$ is uniformly bounded in the space of Radon measures on $[- 1, 1]$. Note that for each $g \in \mathbb{F}_{\infty}$, $\frac{1 - e^{- t \psi(g)}}{t} \to \psi(g)$ as $t \to 0$. Next, we focus on the terms $\nu_{t}$.

Consider $\Lambda = \{ t \in \mathbb{R} : 0 < t \leq 1 \}$ as a directed set with partial order, $\dot{\leq}$ defined as follows: $s \dot{\leq} t$ means that $s > t$. So, $\nu_{t}$ is a net in $M^{+} ([- 1, 1])$, the space of positive Radon measure on $[- 1, 1]$.

Since $(\nu_{t})_{t \in \Lambda}$ is a bounded set in $M([ - 1, 1])$, the space of Radon measures on $[- 1, 1]$, by the Banach-Alaoglu theorem, there exists a subnet $(\nu_{t_{\alpha}})_{\alpha \in E}$ and $\nu \in M([- 1, 1])$ such that $\nu_{t_{\alpha}} \to \nu$ in the weak-* topology of $M([- 1, 1])$.

\begin{align*}
&\int_{- 1}^{1}  \frac{1 - s^{\left\| g \right\|_{2}^{2}}}{1 - s} \ d \nu(s) = \lim_{\alpha} \int_{- 1}^{1}  \frac{1 - s^{\left\| g \right\|_{2}^{2}}}{1 - s} \ d \nu_{t_{\alpha}}(s) \\
&= \lim_{\alpha} \frac{1 - e^{- t_{\alpha} \psi(g)}}{t_{\alpha}} = \lim_{t \to 0} \frac{1 - e^{- t \psi(g)}}{t} = \psi (g)
\end{align*}

\begin{displaymath}
\text{The equality } \lim_{\alpha} \frac{1 - e^{- t_{\alpha} \psi(g)}}{t_{\alpha}} = \lim_{t \to 0} \frac{1 - e^{- t \psi(g)}}{t} \text{ holds true due to the following:}
\end{displaymath}

\begin{itemize}
    \item $t_{\alpha} \dot{\leq} t_{\beta}$ whenever $\alpha < \beta$ (in $E$)
    
    \item For each $r \in \Lambda$, there exits $\alpha \in E$ such that $r \dot{\leq} t_{\alpha}$.
\end{itemize} 

For the direction (2) $\Longrightarrow (1)$, let $\nu$ be a probability measure on $[- 1, 1]$. Note that for each $s \in [- 1, 1]$, the function $g \mapsto \frac{1 - s^{\| g\|_{2}^{2}}}{1 - s}$ is conditionally negative definite. Taking finite sums and limits, we deduce that $g \mapsto \int_{- 1}^{1} \frac{1 - s^{\left\| g \right\|_{2}^{2}}}{1 - s} \ d \nu(s)$ is conditionally negative definite.
\end{proof}

\subsection{Case of $\ell^{p}$ length of the free real line with infinite generators for $0 < p \leq 2$} \

Now, we focus on the case where the group  $G$ is the free real line $\mathbb{R}_{\infty}$ which is defined as $\mathbb{R}_{\infty}:=\Large{*}_{i=1}^{\infty} \mathbb{R}$ where $\Large{*} $ denotes the free product of group and the function $\theta = \left\| \cdot \right\|_{p}^{p}$ is the $\ell^{p}$ length of an element in $\mathbb{R}_{\infty}$ for  $ p\in (1, 2]$.

In this case, $ran(\theta) = \mathbb{R}_{+}$ and hence, $\widehat{ran(\theta)} \cong [0, \infty]$. More precisely, given $\psi \in \widehat{ \mathbb{R}_{+}}$, either there exists a unique $a \in [ 0, \infty )$ such that $\psi(s) := \rho_{a} (s) = e^{- a s}$ or $\psi (s) = \rho_{\infty} (s) := \chi_{\{ 0 \} } (s)$.
We have the following characterization of the positive definite functions and the conditionally negative definite functions on $\mathbb{R}_{+}$, as given in \cite[Proposition 4.4.2, 4.4.3]{BCR1984}.

\begin{lemma}\thlabel{real} 
A function $\varphi : \mathbb{R}_{+} \to \mathbb{R}$ is positive definite and bounded if and only if it has the form

\begin{align*}
    \varphi(s) = \int_{0}^{\infty} e^{- a s} d \mu (a) + b \chi_{\{ 0 \}} (s), s \geq 0,
\end{align*}

where $\mu \in M_{+}^{b} (\mathbb{R}_{+})$ is a bounded positive Radon measure and $b \geq 0$. The pair $(\mu, b)$ is uniquely determined by $\varphi$.

\end{lemma}

\begin{lemma} \thlabel{real negative}
Let $\psi : \mathbb{R}_{+} \to \mathbb{R}$ be a function. Then, $\psi$ is conditionally negative definite and bounded below if and only if it has the form

\begin{displaymath}
\psi(s) = \psi(0) + c s + b \chi_{(0, \infty)} (s) + \int_{0}^{\infty} \left( 1 - e^{- a s} \right) d \mu (a) , \ s \geq 0,
\end{displaymath}

where $b$, $c \geq 0$ and $\mu$, a positive Radon measure on $(0, \infty)$ (possibly infinite), are uniquely determined by $\psi$.
\end{lemma}
Next, we show that the $\ell^{p}$ length function is indeed a partial morphism on $\mathbb{R}_{\infty}$. 

\begin{proposition} \thlabel{realcase} The function $\| \cdot \|_{p}^{p} : \mathbb{R}_{\infty} \to \mathbb{R}_{+}$ is a partial morphism on $\mathbb{R}_{\infty}$.
\end{proposition}
\begin{proof}
Let $M \in \mathbb{N}$, $r_{1}, \cdots , r_{M} \in S$. First, we enumerate the generators as $\{ g_{1, 1}, g_{1, 2}, \cdots , g_{2,1}, g_{2, 2} , \cdots, g_{M, 1}, g_{M, 2}, \cdots \}$.

Now, let $1 \leq j \leq M$. There exists $\lambda_{j} \in \mathbb{R}$ such that $( \lambda_{j} )^{p} = r_{j}$. Then, define $q_{n}: \{ r_{1}, \cdots, r_{M} \} \rightarrow \mathbb{R}_{\infty}$ by $q_{n} (r_{j}) := g_{j, n}^{\lambda_{j}}$. We note the following observation: For each $n, m \in \mathbb{N}$ and $1 \leq j , k \leq M$,

\begin{align*}
\|[q_{n} (r_{j})]^{- 1} q_{m} (r_{k}) \|_{p}^{p} &=
\begin{cases}
0
& \text{if } ~m=n~and~r_{j} = r_{k},\\
r_{j} + r_{k} & \text{otherwise}.
\end{cases}
\end{align*}
\end{proof}

Now, we provide a proof that for each $t > 0$, the function $r \mapsto e^{- t \| r \|_{p}^{p}}$ is positive definite on $\mathbb{R}_{\infty}$.

\begin{proposition}[\cite{Bo86}, Corollary 1]\label{generallength}
Let ${\Bbb F}_q$ (resp. $\mathbb{\mathbb{R}}_{q}$) be the free group (resp. free real line) with generators $r_1, r_2, \cdots, r_{q}$, where $q \in \mathbb{N} \cup \{ \infty \}$. Let $0 < p \leq 2$. Then, for all $t > 0$, the function $\varphi $ defined by $\varphi(x) := e^{- t \|x\|_{p}^{p}}$ is positive definite on ${\Bbb F}_q$ (resp. $\mathbb{\mathbb{R}}_{q}$).
\end{proposition}

\begin{proof} First, let $t > 0$. Observe that the function $\phi : \mathbb{R} \to \mathbb{R}$ defined by $\phi(s) := e^{- t |s|^{p}}$ is positive definite for each $0 < p \leq 2$. The $p=2$ case is well known. The $p<2$ cases follow from the fact that $e^{-|s|^p}$ is an average of $e^{-t|s|^2}$ in $t$. Since $\mathbb{R}_{q} = {\Large{*}}_{i=1}^{q} \mathbb{R}$ (where $\Large{*}$ denotes the free product of groups) and $\varphi(r) = e^{- t \left\| r \right\|_{p}^{p}} = e^{- t |r_{j_{1}}|^{p}} \cdots e^{- t |r_{j_{n}}|^{p}} = \left( \Large{*}_{i = 1}^{n} \phi \right) (r)$ , we deduce that $\varphi$ is a positive definite function. Here we use the fact that a free product of unital positive definite functions is positive definite (see \cite[Corollary 1]{Bo86}).  \end{proof}

Next, we can verify that the function $\psi_{\infty} : \mathbb{R}_{\infty} \to \mathbb{R}_{+}$ defined by: $\psi_{\infty} (r) :=~ \chi_{\{ e \}} (r)$ is positive definite. With all these, we obtain the following characterization.

\begin{corollary} \thlabel{positive free real} Let $0 < p \leq 2$. Given $\varphi : \mathbb{R}_{\infty} \to \mathbb{R}$ an $\ell^{p}$ radial function, the following are equivalent. 

\begin{enumerate}
\item $\varphi$ is positive definite on $\mathbb{R}_{\infty}$.

\item There exist a bounded, positive, Radon measure $\mu$ on $[0, \infty)$ and $b \geq 0$ such that 

\begin{displaymath}
\varphi(r)=\dot{\varphi}(\|r\|_{p}^{p})= \int_{0}^{\infty} e^{- t {\|r\|_{p}^{p}}} \ d \mu(t) + b \chi_{ \{ e \} } (r).
\end{displaymath}

\end{enumerate}

Moreover, if (2) holds, then $\mu$ is uniquely determined by $\varphi$.

\end{corollary}

\begin{proof} (1) $\Longrightarrow$ (2) follows from \thref{main}, \thref{real} and \thref{realcase}. To prove (2) $\Longrightarrow$ (1), let $\mu$ be a bounded, positive, Radon measure on $[ 0, \infty)$. By \thref{positive free real}, for each $t \in [0, \infty)$, the function $\psi_{t} (r) := e^{- t {\| r \|_{p}^{p}}}$ is positive definite on $\mathbb{R}_{\infty}$. Also, the function  defined by $\psi_{\infty} (r) := \chi_{e} (r)$ is positive definite on $\mathbb{R}_{\infty}$. Taking finite sums and limits, we deduce that the function $\varphi : \mathbb{R}_{\infty} \to \mathbb{R}$ defined by $\varphi(r) := \int_{0}^{\infty} e^{- t \| r \|_{p}^{p}}  \ d \mu (t) + b \chi_{e} (r)$ is positive definite on $\mathbb{R}_{\infty}$. \end{proof}

\begin{corollary} Let $0 < p \leq 2$. Let $\psi : \mathbb{R}_{\infty} \to \mathbb{R}$ be an $\ell^{p}$-radial function. Then, the following are equivalent: 

\begin{enumerate}
    \item $\psi$ is conditionally negative definite and bounded below on $\mathbb{R}_{\infty}$.
    
    \item There exist unique $b$, $c \geq 0$ and a positive Radon measure $\nu$ on $(0, \infty)$ (possibly infinite) such that 
    
    \begin{displaymath}
    \psi (r) = \dot{\psi} \left(\| r \|_{p}^{p} \right) = \psi (e) + c \| r \|_{p}^{p} + b \chi_{\mathbb{R}_{\infty} \setminus \{ e \}} (r) + \int_{0}^{\infty} 1 - e^{- t \| r\|_{p}^{p}} \ d \nu (t). 
    \end{displaymath}
\end{enumerate}

\end{corollary}

\begin{proof} $(1 \Longrightarrow 2)$ By \thref{realcase}, the function $\| \cdot \|_{p}^{p} : \mathbb{R}_{\infty} \to \mathbb{R}_{+}$ is a partial morphism on $\mathbb{R}_{\infty}$. By \thref{exact negative definite}, the function $\dot{\psi} : \mathbb{R}_{+} \to \mathbb{R}$ defined by $\dot{\psi} \left( \| r \|_{p}^{p} \right) := \psi (r)$ is conditionally negative definite in the semi-group sense and bounded below on $\mathbb{R}_{+}$. By \thref{real negative}, we obtain (2).

To prove $(2 \Longrightarrow 1)$, we note that the function $\| \cdot \|_{p}^{p} : \mathbb{R}_{\infty} \to \mathbb{R}_{+}$ is conditionally negative definite for all $0 < p \leq 2$ by Schoenberg's theorem and Proposition \ref{generallength}. Also, since the function $\chi_{ \{ e \} } : \mathbb{R}_{\infty} \to \mathbb{R}_{+}$ is positive definite and $\chi_{\mathbb{R}_{\infty} \setminus \{ e \} } = 1 - \chi_{\{ e \}}$, $\chi_{\mathbb{R}_{\infty} \setminus \{ e \} }$ is conditionally negative definite on $\mathbb{R}_{\infty}$. Finally, since $e^{- t \| \cdot \|_{p}^{p}}$ is positive definite for all $t > 0$ by \thref{realcase}, $1 - e^{- t \| \cdot \|_{p}^{p}}$ is conditionally negative definite on $\mathbb{R}_{\infty}$ for all $t > 0$. Taking finite sums and limits, $\int_{0}^{\infty} 1 - e^{- t \| \cdot \|_{p}^{p}} \ d \nu (t)$ is conditionally negative definite on $\mathbb{R}_{\infty}$.
\end{proof}

\subsection{Case of $\ell^{p}$ length of $\mathbb{R}^{\mathbb{N}}$ for $0 < p \leq 2$} \

Now, we focus on the case of the group $G = \mathbb{R}^{\mathbb{N}}$, the infinite direct product of countably many copies of $\mathbb{R}$ and the function $\theta = \left\| \cdot \right\|_{p}^{p}$ is the $\ell^{p}$ length of an element in $\mathbb{R}^{\mathbb{N}}$, where $0 < p \leq 2$. More precisely, $\mathbb{R}^{\mathbb{N}} := \{ (a_{n})_{n = 1}^{\infty} : (a_{n})_{n = 1}^{\infty}$ has finite support$\}$.  The proof is essentially similar to the case for the free real line with infinite generators. We will only state the theorems whose proofs are similar to the previous case.

\begin{proposition} \thlabel{realcasecommutative} The function $\| \cdot \|_{p}^{p} : \mathbb{R}^{\mathbb{N}} \to \mathbb{R}_{+}$ is a partial morphism $\mathbb{R}^{\mathbb{N}}$.
\end{proposition}

\begin{corollary} \thlabel{positive equivalence real} Let $0 < p \leq 2$. Given $\varphi : \mathbb{R}^{\mathbb{N}} \to \mathbb{R}$ an $\ell^{p}$ radial function, the following are equivalent. 

\begin{enumerate}
\item $\varphi$ is positive definite on $\mathbb{R}^{\mathbb{N}}$.

\item There exist a bounded, positive, Radon measure $\mu$ on $[0, \infty)$ and $b \geq 0$ such that 

\begin{displaymath}
\varphi(r)=\dot{\varphi}(\|r\|_{p}^{p})= \int_{0}^{\infty} e^{- t {\|r\|_{p}^{p}}} \ d \mu(t) + b \chi_{ \{ e \} } (r).
\end{displaymath}

\end{enumerate}

Moreover, if (2) holds, then $\mu$ is uniquely determined by $\varphi$.

\end{corollary}

\begin{corollary} \thlabel{negative equivalence real} Let $0 < p \leq 2$. Let $\psi : \mathbb{R}^{\mathbb{N}} \to \mathbb{R}$ be an $\ell^{p}$-radial function. Then, the following are equivalent: 

\begin{enumerate}
    \item $\psi$ is conditionally negative definite and bounded below on $\mathbb{R}^{\mathbb{N}}$.
    
    \item There exist unique $b$, $c \geq 0$ and a positive Radon measure $\nu$ on $[0, \infty)$ (possibly infinite) such that 
    
    \begin{displaymath}
    \psi (r) = \dot{\psi} \left(\| r \|_{p}^{p} \right) = \psi (0) + c \| r \|_{p}^{p} + b \chi_{\mathbb{R}^{\mathbb{N}} \setminus \{ 0 \}} (r) + \int_{0}^{\infty} \left( 1 - e^{- t \| r\|_{p}^{p}} \right) d \nu (t) .
    \end{displaymath}
\end{enumerate}

\end{corollary}

With these results, we obtain the classical Schoenberg-Bochner theorem (see e.g. \cite[Theorem 13.14]{SSV2012}) as a corollary.

\begin{corollary} \thlabel{commutative result}
Let $f : [0, \infty) \to [0, \infty)$ be a function. Then, the following are equivalent.

\begin{enumerate}
    \item There exists $0 < p \leq 2$ such that for all $d \in \mathbb{N}$, the function $\xi \mapsto f(|\xi|^{p})$, $\xi \in \mathbb{R}^{d}$ is conditionally negative definite.  
    
    \item f is conditionally negative definite in the semi-group sense on $[0, \infty)$.
    
    
    
\end{enumerate}   
\end{corollary}

\begin{proof} (1) $\Longrightarrow$ (2). If $f \circ \| \cdot \|_{p}^{p} : \mathbb{R}^{d} \to \mathbb{R}$ is conditionally negative definite and bounded below for all $d \in \mathbb{N}$, then $f \circ \| \cdot \|_{p}^{p} : \mathbb{R}^{\mathbb{N}} \to [ 0, \infty)$ is conditionally negative definite and bounded below on $\mathbb{R}^{\mathbb{N}}$ because the sum in Definition 2.8 is a finite sum. By \thref{exact} and \thref{realcase}, (2) is true.

(2) $\Longrightarrow$ (3) is \thref{real negative}.




\end{proof}

\section*{Acknowledgement}
The authors would like to thank Yazhou Han for his advice and discussion for the project.  The authors are partially supported by the NSF Grant DMS 1700171.

\bibliographystyle{amsplain}

\end{document}